\input amstex
\documentstyle{amsppt}
\magnification=\magstep1
\hsize=5in
\vsize=7.3in
\TagsOnRight
\topmatter
\title Applications of  braided endomorphisms from conformal inclusions
\endtitle
\author Feng Xu \endauthor
\abstract We give three applications of general theory about
braided endomorphisms  from conformal inclusions developed
previously by us. The first is an example of subfactors associated
with conformal inclusion whose dual fusion ring is non-commutative.
In the second application we show that the Kac-Wakimoto hypothesis
about certain relations between  branching rules and S-matrices,
which has existed for almost a decade, is {\bf not} true in
at least three examples. Finally we show that the fusion rings
of subfactors associated with the conformal inclusions 
$SU(n)_{(n+2)}\subset SU(n(n+1)/2)$ and  $SU(n+2)_n\subset SU((n+1)(n+2)/2)$
are canonically isomorphic using a version of level-rank duality.
\endabstract
 
\thanks This work originated from a very informative discussion
with Prof.Victor Kac at the International Congress of Mathematical
physics 97 in Brisbane, Australia. I'd like to thank Prof.Victor Kac
for pointing out the reference [LL] and for bringing my attention to
Kac-Wakimoto hypothesis.  This work is partially
supported by NSF grant DMS-9500882.\endthanks
\endtopmatter
\heading \S0.  Introduction \endheading
Let $G_k \subset H$ be a conformal inclusion with both $G$ and $H$
being semisimple compact Lie groups, and $k$ the Dykin index of the
inclusion (cf. [GNO] or
[KW]). We shall use $i$ (resp. $\lambda$) to denote the irreducible
projective positive energy representation of loop group $LH$ (resp.
$LG$) at level 1 (resp. k) (cf.[PS]).  Denote by $b_{i \lambda}$ the branching
coefficients, i.e., when restricting to $LG$, $i$ decomposes as
$\sum_{\lambda} b_{i \lambda} \lambda$.  Denote by $S_{ij}$ (resp. 
$S_{\lambda \mu}$) the S-matrices of  $LH$ (resp.
$LG$) at level 1 (resp. k) (cf. [Kac2]).  A general hypothesis of 
 Kac-Wakimoto (cf. [KW]) in this case states that
$S_{\lambda \mu} \bar S_{ij} \geq 0$ whenever $b_{i \lambda} b_{j \mu}\neq 0$.
This hypothesis has been checked to be correct in many cases.\par
Although there are explicit formula for $S$ matrices (cf. [Kac2]), 
it is in general
very difficult to calculate  $S$ matrices in the high level case. 
For this reason there is no general theoretical approach to  Kac-Wakimoto
hypothesis except case by case study. \par
In this paper we will argue that the general theory developed in [X1],
which is motivated by subfactor theory, is particularly useful in 
studying the `` Kac-Wakimoto hypothesis type" question.  We find three
counter examples to  Kac-Wakimoto hypothesis, and in the course of
studying such questions, we also find answers to some of the questions
naturally arised in [X1] and [X3].\par
Let us describe the content of this paper in more details. \par
In \S 2 we recall some of the material in [X1] to set up notations.
In \S2.1 we define sectors and what we mean by fusion ring and
dual fusion ring associated with a sector. In \S2.2 we recall the 
braided endomorphisms from conformal inclusions and summarize their
properties in Th.2.2 and Th.2.3.  Both theorems are proved in [X1]. \par
\S 3 and \S4 contain the main results of this paper which are proved
by using the results stated in \S 2.

In \S3, we first notice a curious inequality
as the first application of results in [X1]. This
inequality can be stated entirely in terms of S-matrices and branching
coefficients, and it seems to be hard to prove it by other means.
 Prop.3.1 is a slightly generalized version of (1) of Th.3.9 in [X1]
about the spectrum of certain ring which arises naturally in  [X1].
 Cor.3.2 shows that the framework of [X1] is almost tailor made
to attack the  Kac-Wakimoto hypothesis.  Prop.3.3 allows one to
obtain information about the non-commutativity of certain elements
in a ring defined in [X1] from the branching rules. \par
Th.3.4 is the main result of \S3 which follows from
applying Cor.3.2 and Prop.3.3 to a special case of conformal
inclusion $SU(3)_9 \subset E_6$.  
In (1) of Th.3.4 we show that certain ring $A_{\rho}$ which appears
natural in [X1] is not generated by the descendants of the
braided endomorphisms.  \par
 In (2) of Th.3.4 we show that the dual
fusion ring of the Jones-Wassermann subfactor associated with
 $SU(3)_9 \subset E_6$ is non-commutative by exhibiting explicitly
non-commutative relations.  Notice it is in general very difficult
to obtain detailed information about 
 the dual fusion ring and this is the first example of
conformal inclusions which
shows non-commutative dual fusion ring. \par
In (3) of Th.3.4 we show the same 
conformal inclusion gives a counter example to the  Kac-Wakimoto hypothesis.
After discovering the counter example in  (3) of Th.3.4 we find two 
more counter examples where $S$ matrices can be computed explicitly by
elementary computations, and they are given at the end of \S3. \par
In \S4.1 we give a version of level-rank duality which follows
naturally from [X1], and this is summarized in Th.4.1. The coincidence
of certain fusion coefficients in the fusion rings associated with
$SU(m)_n$ and $SU(n)_m$ are already noticed before, but Th.4.1 gives
an embedding of these two fusion rings in one ring where they can
be compared.  By using Th.4.1, we show in Th.4.2 that the fusion rings
of the Jones-Wassermann  subfactors associated with the conformal inclusions
$SU(n)_{(n+2)}\subset SU(n(n+1)/2)$ and  $SU(n+2)_n\subset SU((n+1)(n+2)/2)$
are canonically isomorphic.  This fact was first noticed in an example
in [X3].\par
In \S5, we give our conclusions and questions. \par
\subheading {2.1.  Sectors} 
Let $M$ be a properly infinite factor
and  $\text{\rm End}(M)$ the semigroup of 
the unit preserving endomorphism of $M$.  In this paper $M$ will always
be the unique hyperfinite $III_1$ factors. 
Let $\text{\rm Sect}(M)$ denote the quotient of $\text{\rm End}(M)$ modulo 
unitary equivalence in $M$. We shall denote by $[\rho]$ the image of
$\rho \in \text{\rm End}(M)$ in  $\text{\rm Sect}(M)$.\par
 It follows from
\cite{L3} and \cite{L4} that $\text{\rm Sect}(M)$, with $M$ a properly
infinite  von Neumann algebra, is endowed
with a natural involution $\theta \rightarrow \bar \theta $ that commutes
with all natural operations ;  moreover,  $\text{\rm Sect}(M)$ is
naturally a semiring.\par
Suppose $\rho \in \text{\rm End}(M)$ is given together with a normal
faithful conditional expectation
$\epsilon:
M\rightarrow \rho(M)$.  We define a number $d_\epsilon$ (possibly
$\infty$) such that (cf. [PP]):
$$
d_\epsilon^{-2} :=\text{\rm Max} \{ \lambda \in [0, +\infty)|
\epsilon (m_+) \geq \lambda m_+, \forall m_+ \in M_+
\}$$.\par
Now assume  $\rho \in \text{\rm End}(M)$ is given together with a normal
faithful conditional expectation
$\epsilon:
M\rightarrow \rho(M)$, and assume  $d_\epsilon < +\infty$. We define
$$
d = \text{\rm Min}_\epsilon \{ d_\epsilon \}
$$.   $d$ is called the statistical dimension of  $\rho$. It is clear
from the definition that  the statistical dimension  of  $\rho$ depends only
on the unitary equivalence classes  of  $\rho$. 
The properties of the statistical dimension can be found in
[L1], [L3] and  [L4].\par
Denote by $\text{\rm Sect}_0(M)$ those elements of 
$\text{\rm Sect}(M)$ with finite statistical dimensions.
For $\lambda $, $\mu \in \text{\rm Sect}_0(M)$, let
$\text{\rm Hom}(\lambda , \mu )$ denote the space of intertwiners from 
$\lambda $ to $\mu $, i.e. $a\in \text{\rm Hom}(\lambda , \mu )$ iff
$a \lambda (x) = \mu (x) a $ for any $x \in M$.
$\text{\rm Hom}(\lambda , \mu )$  is a finite dimensional vector 
space and we use $\langle  \lambda , \mu \rangle$ to denote
the dimension of this space.  $\langle  \lambda , \mu \rangle$
depends 
only on $[\lambda ]$ and $[\mu ]$. Moreover we have 
$\langle \nu \lambda , \mu \rangle = 
\langle \lambda , \bar \nu \mu \rangle $, 
$\langle \nu \lambda , \mu \rangle 
= \langle \lambda , \mu \lambda \rangle $ which follows from Frobenius 
duality (See \cite{L2} or \cite{Y}).  We will also use the following 
notations: if $\mu $ is a subsector of $\lambda $, we will write as
$\mu \prec \lambda $  or $\lambda \succ \mu $.  A sector
is said to be irreducible if it has only one subsector. \par
For an endomorphism $\rho \in \text{\rm End}(M)$, the fusion ring
(resp. the dual fusion ring) associated with the inclusion 
$\rho(M) \subset M$ is defined to be the ring generated by the
irreducible descendants of $\rho \bar \rho$ (resp. $\bar \rho \rho$).
The origin of such notions comes from subfactor theory (cf.[DB]).
\vskip .1in
\noindent
\subheading {2.2.  Braided endomorphisms from conformal inclusions} 
In this paper, we shall restrict our attention to the following
conformal inclusions (cf. [PZ] or [GNO]):
 $ \ \ {SU}(2)_{10} \subset {Spin}(5), \
{SU}(2)_{28} \subset  G_2$, 
 
\noindent
 $ \ \ {SU}(3)_5 \subset {SU}(6), \
{SU}(3)_9 \subset  E_6, \ {SU}(3)_{21} \subset
 E_7$; \par
\noindent
$( A_8)_1 \subset ( E_8)$;\par
\noindent   
and five infinite series:
$$
\align
{SU}(N)_{N-2} & \subset \ {SU} \left( \frac{N(N-1)}{2}
\right), \ \ N \geq 4 ; \tag a \\
{SU}(N)_{N+2} & \subset \ {SU} \left( \frac{N(N+1)}{2}
\right),  \tag b \\
{SU}(N)_{2N} & \subset \ {Spin} (4N^2 - 1), \ \ N \geq
2;  \tag c \\
{SU}(2N+1)_{2N+1} & \subset \ {Spin}(4N(N+1)); \tag d \\
{SU}(M)_N \times {SU}(N)_M & \subset \ {SU}(NM).
\tag e\endalign
$$
These cover all the maximal conformal inclusions of the
form $G=SU(N) \times SU(M) \subset H$ with
$H$ being a simple simply connected group.\par
The aim of this section is to recall some of the results in [X1] which
will be used in \S3 and \S4. For the proofs and unexplained terminology,
we safely refer the reader to [X1].
For the representation theory of Loop groups, we refer the reader
to [PS].
\par
  Denote by  $\lambda$                  
a positive energy projective representation of $LG$ at the level (=
Dykin index given above).   $\lambda$ is not necessarily irreducible but
we assume that it is a finite direct sum of irreducible representations.
By [W2], each $\lambda$  naturally gives a sector, denoted by the same
letter  $\lambda \in \text{\rm Sect}(M)$, where $M$ is the 
unique hyperfinite $III_1$ factor.  Such sectors generate a finite
dimensional ring denoted by $Gr(C_k)$ with $k$ indicating the level. 
 A basis of $Gr(C_k)$ is given by all the irreducible 
 positive energy projective representation, denoted by
 $\mu_i$'s, of $LG$ of fixed level.  The
structure constants $N_{\mu_1\mu_2}^{\mu_3}$ are given by:
$$
\mu_1\mu_2 = \sum_{\mu_3} N_{\mu_1\mu_2}^{\mu_3} \mu_3
$$, and it is known (cf. Cor.1 in Chapter 34 of [W2] and P.288 of [Kac2])
that  $N_{\mu_1\mu_2}^{\mu_3}$ are determined uniquely by S-matrices
of $LG$. \par
Let $\gamma_i:= \sum_i b_{i \lambda} \lambda$. We shall use $1$ to
denote the vacuum representation of $LH$. It is shown (cf. (1) of
Prop.2.8 in [X1]) that there are sectors $\rho, \sigma_i \in 
\text{\rm Sect}(M)$ such that:
$$
\rho \sigma_i \bar \rho =\gamma_i
$$.  
Notice that  $\sigma_i$ are in one-to-one correspondence with
the irreducible representations $i$ of $LH$ and they are 
irreducible sectors, generating a finite dimensional ring.
 The structure constants $N_{ij}^k$ are given by :
$$
\sigma_i\sigma_j = \sum_k N_{ij}^k \sigma_k
$$,  and  $N_{ij}^k$ are uniquely determined 
the S matrices of $LH$ (cf. the paragraph after Th.1.6. in [X1]). \par
The subfactor $\rho(M) \subset M$ is called the Jones-Wassermann 
subfactor associated with conformal inclusions.  For more detailed
discussions on  the Jones-Wassermann
subfactor, we refer the reader to [X1]. \par
The crucial observation in [X1] is the following:
for each $\lambda$, there exists a sector denoted 
\footnotemark\footnotetext{In
the notations of [X1] it should be denoted by $[a_\lambda]$ to 
emphasize that it is a sector rather than an endomorphism, but we
omit $[,]$ for simplicity.}  by $a_\lambda$ , 
such that the following theorem is true (cf. Th.3.1, Cor.3.2 and
Th.3.3 of [X1]):
\proclaim{Theorem 2.2}
(1).
The map $\lambda \rightarrow a_\lambda$ is a ring homomorphism;\par
(2). 
$\rho a_\lambda = \lambda \rho, a_\lambda \bar \rho = \bar \rho \lambda$;\par
(3). $\langle \rho a_\lambda,  \rho a_\mu \rangle =
\langle a_\lambda,  a_\mu \rangle = \langle  a_\lambda \bar \rho,  a_\mu
\bar \rho  \rangle$; \par
(4).  $\langle \rho a_\lambda, \rho \sigma_i  \rangle =
\langle a_\lambda, \sigma_i  \rangle = 
\langle a_\lambda \bar \rho, \sigma_i \bar \rho
\rangle $ ; \par
(5). (3) (resp. (4)) remains valid if $ a_\lambda,  a_\mu$ (resp.  $a_\lambda$)
is replaced by any of its subsectors.
\endproclaim
\demo{Proof}
(1) to (4) follows from  Th.3.1, Cor.3.2 and
Th.3.3 of [X1], so we just have to show (5).  Let $b$ (resp. $c$) be
subsectors of  $ a_\lambda$ (resp. $ a_\mu$), then we have:
$$
\langle \rho b, \rho c \rangle \geq \langle b,  c \rangle
$$.  But if:
$$
\langle \rho b, \rho c \rangle > \langle b,  c \rangle
$$, we would have 
$$
\langle \rho a_\lambda,  \rho a_\mu \rangle > \langle a_\lambda,  
a_\mu \rangle
$$, contradicting (3).  The remaining cases are proved similarly.
\enddemo
\hfill Q.E.D.
\par 
Perhaps the most surprising part is (3) and (4) of Th.2.2. Notice one obviously 
has $ \langle \rho a_\lambda,  \rho a_\mu \rangle \geq
\langle a_\lambda,  a_\mu \rangle$, the nontrivial part is the equality
which is proved in [X1] by locality considerations of [LR].  Combined
with Frobenius duality, Th.3.3 puts a powerful constraint on the
subsectors of $ a_\lambda$.\par
To make a connection with the branching coefficients introduced in \S1,
we have:
$$
\align
\langle a_\lambda, \sigma_j \rangle & = \langle a_\lambda \bar \rho, 
\sigma_j \bar \rho \rangle  \\
& = \langle \bar \rho \lambda, \sigma_j \bar \rho \rangle \\
& =  \langle \lambda, \rho \sigma_j \bar \rho \rangle \\
& = \langle \lambda, \gamma_j \rangle = b_{j\lambda}
\endalign
$$, where we have used (3) and (2) of Th.2.2 in the first
and the second step, Frobenius duality in the third step, 
the equation $ \rho \sigma_j \bar \rho  = \gamma_j$ in the
fourth step and definitions
of $\gamma_j$ in the last step. \par
It is shown (cf. \S3.4 [X1]) there is another sector $\tilde  a_\lambda$
with exactly the same properties as $ a_\lambda$ in Th.2.4. It is 
associated to the choice of under-crossing of braiding in the definition
 and it is known that  $\tilde  a_\lambda$ is in general different from
 $ a_\lambda$ (cf. Lemma 3.2 of [X1]).\par
Define $C_\rho$ (resp. $\tilde C_\rho$) to be the complex 
finite dimensional ring
generated by the irreducible subsectors of  $ a_\lambda$ (resp. 
 $\tilde  a_\lambda$) for all $\lambda$ of fixed level.  The paring
$\langle, \rangle$ introduced in \S 2.1 extends by linearily in
the first variable and conjugate linearily in the second variable
to a positive definite form on  $C_\rho$. Notice the conjugation
$b\rightarrow \bar b$  in \S 2.1  extends   conjugate linearily to
the elements in $C_\rho$ such  that the Frobenius duality holds.  \par
Define
$A_\rho$  to be the finite dimensional ring
generated by the irreducible subsectors of $\bar \rho \lambda \rho$
 for all $\lambda$ of fixed level. Notice that
$C_\rho \subset A_\rho$.  The commutativity of certain elements
in these rings are investigated and they are summarized in the following
theorem (cf. Th.3.6, Lemma 3.3 of [X1])
\proclaim{Theorem 2.3}
 Let $b$ be any subsector of $A(a_\mu)$
where $A$ is an arbitrary polynomial in $a_\mu, \ \mu \in Gr(C_k)$, then
$a_\lambda b = ba_\lambda$ for any $\lambda \in  Gr(C_k)$;
 
(2) Let $c$ be any subsector of $\bar\rho \rho$, then $a_\lambda c
= ca_\lambda$ for any $\lambda \in  Gr(C_k)$; \par
(3)  Let $x$ and $y$ be subsectors of
$\tilde a_\lambda$
and $a_{\mu}$ respectively.  Then $xy = yx$.
\endproclaim
\heading \S3.  An example of conformal inclusion \endheading
We preserve the set up of \S2.2. 
We shall denote the set of irreducible sectors of $a_\lambda$ for
all $\lambda$ of fixed level by $V$.  Notice $\sigma_i \in V$,
and these are refered to as "special nodes" in \S3.4 of [X1].  Let:
$$
 a_\lambda  a
= \sum_{b\in V} V^\lambda_{ab} b
$$, where $V^\lambda_{ab}$ are nonnegative
integers.  Denote by $V^\lambda$ the matrix such that $(V^\lambda)_a^b = 
V^\lambda_{ab}$.  By (1) of Th.2.2 $V^{\mu_1} V^{\mu_2} = 
\sum_{\mu_3} N_{\mu_1\mu_2}^{\mu_3} V^{\mu_3},$  so we have:
$$
\sum_{\mu_3} N_{\mu_1\mu_2}^{\mu_3} V^{\mu_3}_{1\sigma_i} =
\sum_a V^{\mu_1}_{1a}  V^{\mu_2}_{a\sigma_i} \geq
\sum_j  V^{\mu_1}_{1 \sigma_j}  V^{\mu_2}_{ \sigma_j\sigma_i} 
$$, where $1$ denotes the identity sector.  Recall that 
$ V^{\mu_3}_{1 \sigma_i} = b_{i \mu_3},  V^{\mu_1}_{1 \sigma_j} 
=  b_{j \mu_1}$,  and 
$$
\align
 V^{\mu_2}_{ \sigma_j\sigma_i} & = 
\langle a_{\mu_2} \sigma_j, \sigma_i \rangle = 
\langle a_{\mu_2} , \sigma_{\bar j} \sigma_i \rangle \\
& = \langle a_{\mu_2} , \sum_k N_{\bar {j} i}^k \sigma_k \rangle \\
& = \sum_k N_{\bar {j} i}^k  \langle a_{\mu_2} ,\sigma_k \rangle \\
& = \sum_k N_{\bar {j} i}^k  b_{k\mu_2}
\endalign
$$.  Putting the above together,  we have derived the following inequality:
$$
\sum_{\mu_3} N_{\mu_1\mu_2}^{\mu_3}  b_{i\mu_3} \geq
\sum_{k,j} N_{\bar {j} i}^k  b_{k\mu_2} b_{j\mu_1}   
$$.  Since $ N_{\mu_1\mu_2}^{\mu_3}$ (resp.$ N_{\bar {j} i}^k$) are determined
uniquely by S-matrices of $LG$ (resp. $LH$) by so called Verlinde formula,
the above inequality is a relation between branching rules and S-matrices,
and we don't know if one can prove it by the usual representation
theoretical approach, i.e., by using the results of [KW] and [Kac2]. \par
Define matrix $N_c$ by $N_{ca}^b = \langle ca,b\rangle$ for
 $a,b,c \in V$.
Then $V^\lambda = \sum_c V^\lambda_{1c} N_c$.
Since $[a_{\bar\lambda}] = [\bar a_\lambda], [\sigma_j a_\lambda]
= [a_\lambda \sigma_j], \ V^\lambda, \ N_{\sigma_j}$ are
commuting normal matrices, so they can be simultaneously diagonalized.
Recall the irreducible representations of $Gr(C_k)$ are given by
$$
\lambda \rightarrow \frac{S_{\lambda \mu}}{S_{1\mu}}.
$$
Assume $V^\lambda_{ab} = \sum_{\mu,i,s\in(\text{\rm Exp})}
\frac{S_{\lambda \mu}}{S_{1\mu}} \cdot \psi_a^{(\mu,i,s)}
\psi_b^{(\mu,i,s)^*}$ where $ \psi_a^{(\mu,i,s)}$ are normalized orthogonal
eigenvectors of $V^\lambda$ (resp. $N_{\sigma_i}$) with eigenvalue
$\frac{S_{\lambda \mu}}{S_{1\mu}}$ (resp.
 $\frac {S_{i j}}{S_{1 j}}$) .
$(Exp)$ is a set of $\mu,i,s$'s
 and $s$ is an
index indicating the multiplicity of  $\mu,i$.
Recall if a representation is denoted by $1$, it will always be the
vacuum representation.
\proclaim{Proposition 3.1}
$(\delta,k,s)\in (Exp)$ if and only if $b_{k\delta}>0$.  Moreover, there is 
a choice of eigenvectors such that  $ \psi_1^{(\delta,k,s)} > 0$ for any
$(\delta,k,s)\in (Exp)$.
\endproclaim
\noindent
{\it Proof:}
Since $b_{j\lambda} = V^\lambda_{1 \sigma_j}$, we have:
$$
b_{j\lambda} = \sum_{(\mu,i,s)\in (Exp)} \frac{S_{\lambda \mu}}{S_{1\mu}}
\frac{S_{\bar j i}}{S_{1i}} |\psi_1^{(\mu,i,s)}|^2
$$,
where we have also used 
$$
\psi_{\sigma_j}^{(\mu,i,s)} = \frac{S_{ j i}}{S_{1i}}
\psi_1^{(\mu,i,s)}
$$ (cf.(4) of Th.3.9 in [X1]).
By using the following equation (cf. [KW]):
$$
\sum_j S_{kj} b_{j\lambda} = \sum_{\mu} b_{k\mu} S_{\mu \lambda}
$$
we obtain:
$$
\sum_{(\mu,s) \in Exp(k)}  \frac{S_{\lambda \mu}}{S_{1\mu}}
\frac{1}{S_{1k}}  |\psi_1^{(\mu,k,s)}|^2
= \sum_{\mu} b_{k\mu} S_{\mu \lambda} 
$$, where $(\mu,s) \in Exp(k)$  means
$k$  is fixed and $(\mu,k,s) \in (Exp)$ .  
Multiply both sides by $\bar S_{\lambda \delta}$ and summed over 
$\lambda$, using the fact that $S$ matrices are 
unitary (cf. [Kac2]), we get:
$$
\sum_{s\in Exp(\delta,k)} \frac{1}{S_{1\delta}S_{1k}}
 |\psi_1^{(\delta,k,s)}|^2 = b_{k\delta} \tag 1
$$, where $s\in Exp(\delta,k)$ means $\delta,k$ are fixed and
$(\delta,k,s) \in (Exp)$.
It follows immediately that if $ b_{k\delta} >0$, then
$(k,\delta,s)\in (Exp)$ for some $s$.\par
Let $(\delta, k, 1),...(\delta, k, p)$ be the subset in $(Exp)$ with
$(\delta,k)$ fixed.  By Th.2.2 and the
comments after it   $\sigma_i$ appears as subsectors of some
$\tilde a_\mu$, so by Th.2.3   ${\bar a}$ commutes with 
$a_\lambda$ and $\sigma_i$ ,
therefore ${\bar a}$ preserves the subspace spanned by vectors 
 $\psi^{(\delta, k, 1)},...\psi^{(\delta, k, p)}$. 
Note for any $a\in V$, we have:
$$
\align
\psi_a^{(\delta,k,s)} & = \langle \psi^{(\delta,k,s)}, a \rangle 
= \langle N_{\bar a}\psi^{(\delta,k,s)}, 1 \rangle \\
& = \sum_t  N_{\bar a s}^t \psi_1^{(\delta,k,t)}
\endalign
$$
where $ N_{\bar a s}^t = \langle N_{\bar a}\psi^{(\delta,k,s)},
\psi^{(\delta,k,t)}\rangle$. \par
It follows that if $(\delta,k,s)\in (Exp)$, then
$\psi_a^{(\delta,k,s)} \neq 0$ for some $a\in V$ which implies
$\psi_1^{(\delta,k,t)} \neq 0$ for some $t$.  By equation (1)
this implies  $b_{k\delta} > 0$. \par
Let $(\delta, k, 1),...(\delta, k, p)$ be the subset in $(Exp)$ with
$(\delta,k)$ fixed. It follows from (1) that we can always make a gauge
choice such that $\psi_1^{(\delta,k,1)} = ... = \psi_1^{(\delta,k,p)} > 0$.
 \hfill Q.E.D. \par
%\nonident
We remark that the existence of the choice of eigenvectors in Prop.3.1
is postulated as (2) of ix in [PZ]. \par
Recall the Kac-Wakimoto hypothesis as stated in the beginning of
the introduction.  It follows immediately from Prop.3.1. that:
\proclaim{Corollary 3.2}
 Kac-Wakimoto hypothesis is true if and only if for any $\lambda$
with $b_{j\lambda} \neq 0$, $V^\lambda N_{\sigma_{\bar j}}$ is
semi-positive definite, i.e. 
$\langle a_\lambda \sigma_{\bar j} x, x \rangle \geq 0$
for any $x\in C_\rho$.
\endproclaim
Cor.3.2 is very effective in verifying that  Kac-Wakimoto hypothesis
is true in many cases. The strategy is to write 
$a_\lambda \sigma_{\bar j} = \sum_b b \bar b$ by using the information
on the ring structure of $C_\rho$ which can be effectively determined
by Th.2.2 in many  cases (cf. examples in [X1]).  We have done so in 
examples of \S4.1 and example 1 and 4 in [X1], and series (e) of
\S2.2.\par
In trying to prove  Kac-Wakimoto hypothesis by using
Cor.3.2,  we find the following proposition:
\proclaim{Proposition 3.3}
(1) If $\sigma_j \bar \rho \rho = \bar \rho \rho \sigma_j$ 
, then $a_{\gamma_j} = a_{\gamma_1} \sigma_j$; \par
(2) $a_{\gamma_j} = a_{\gamma_1} \sigma_j$ for all $j$ if and only
if the following holds: if for any $i,\mu$ with  $b_{i\mu} >0$, we
have $b_{j\mu} =0$ for $j\neq i$.
\endproclaim
\noindent
{\it Proof:}
(1):
It is sufficient to show that 
$$
\langle a_{\gamma_j}, b \rangle = \langle  a_{\gamma_1} \sigma_j, b \rangle
$$ for any $b\in V$. We have:
$$
\align
\langle a_{\gamma_j}, b \rangle & = \langle a_{\gamma_j} \bar \rho, b \bar \rho
\rangle \\
&= \langle a_{\gamma_j} \bar \rho \rho, b\rangle \\
&= \langle \bar \rho \gamma_j \rho, b\rangle \\
&= \langle \bar \rho \rho \sigma_j \bar \rho \rho, b \rangle
\endalign
$$ where we have used (5) of Th.2.2 in the first step, Frobenius duality in the 
second step and (2) of Th.2.2 in the third step and the identity
at the beginning of \S3 in the last step.\par
On the other hand, we have:
$$
\align
\langle a_\gamma \sigma_j, b \rangle & = \langle  a_\gamma \sigma_j \bar \rho,
b \bar \rho \rangle \\
&= \langle a_\gamma \sigma_j \bar \rho \rho,b \rangle \\
&=  \langle \sigma_j a_\gamma \bar \rho \rho,b \rangle \\
&=  \langle \sigma_j \bar \rho \rho \bar \rho \rho, b \rangle \\
&=  \langle \bar \rho \rho \sigma_j \bar \rho \rho, b \rangle
\endalign
$$
where the intermediate steps are similar to the previous one and
in the last step we have used the condition 
$\sigma_j \bar \rho \rho = \bar \rho \rho\sigma_j$. \par
(2)
Notice that  $a_{\gamma_j} = a_{\gamma_1} \sigma_j$ iff 
$V^{\gamma_j} = V^{\gamma} N_{ \sigma_j}$. 
It follows from Prop.3.1 that:
 $a_{\gamma_j} = a_{\gamma_1} \sigma_j$ for all $j$ iff
for any $i,\mu$ with $b_{i\mu} \neq 0$, we have
$$
\sum_\lambda b_{j\lambda} \frac{S_{\lambda\mu}}{S_{1\mu}}
=\sum_\lambda b_{1\lambda}  \frac{S_{\lambda\mu}}{S_{1\mu}}
\frac{S_{ji}}{S_{1i}} \tag 2
$$.
Use the fact $\sum_\lambda b_{j\lambda} S_{\lambda\mu} =
\sum_k S_{jk} b_{k\mu}$, (2) is equivalent to:
$$
\sum_k  S_{jk} b_{k\mu} = \sum_k S_{1k} b_{k\mu} 
\frac{S_{ji}}{S_{1i}}
$$.
Multiply both sides of the above equation 
by $\bar S_{jl}$ and summed over $j$,
it is easy to see that (2) is equivalent to
the statement that $b_{l\mu} =0$ for any $l\neq i$. \hfill Q.E.D.
Now we are ready to apply Cor.3.2 and Prop.3.3 to an example,
i.e., the conformal inclusion 
 ${SU}(3)_9 \subset {E_6}$ (cf. example 2 in \S4 of [X1]). 
Denote by $H_0, H_1, H_2$ the level 1 irreducible representations
of $LE_6$ ($H_1$ is the vacuum
representation) and label the dominant weights of $SU(3)$ by $(\lambda_1,
\lambda_2)$ with $\lambda_1\geq \lambda_2 \geq 0$, we have the
following decompositions:
$$
\align
H_0 & = H_{(4,2)} + H_{(7,2)} + H_{(7,5)} \\
H_2  & = H_{(4,2)} + H_{(7,2)} + H_{(7,5)} \\
H_1 &= H_{(0,0)} + H_{(9,0)}
+     H_{(9,9)} + H_{(8,4)}+  H_{(5,1)} +  H_{(5,4)}
\endalign
$$.
By using similar calculations as in  example 2 in \S4 of [X1] and 
the computation
of statistical dimensions, one easily obtains the following identities:
$$
\align
a_{(2,1)}^2 & = 6 a_{({2,1})} + \sigma_0 + \sigma_1 +  \sigma_2 \tag 3\\
a_{(4,2)} & = 2  a_{({2,1})} + \sigma_0  +  \sigma_2 \tag 4 \\
a_{(5,1)} & =  2  a_{({2,1})} + \sigma_1 \tag 5
\endalign
$$. In fact, $d_{a_{(2,1)}} = 3 + 2\sqrt 3, d_{a_{(4,2)}} = 8 + 4 \sqrt 3$
and the above identities are proved by first showing that the left side
contains the right side using Th.2.2
and then by comparing the  statistical dimensions.
Now we are ready to prove the following theorem:
\proclaim{Theorem 3.4}
For the  conformal inclusion
 ${SU}(3)_9 \subset {E_6}$, we have: \par
(1) The ring $A_\rho$ is not generated by the subsectors of
$a_\lambda$ and $\tilde a_\lambda$ for any dominant weights $\lambda$ 
of  ${SU}(3)$ at level 9; \par
(2) The dual fusion ring of the Jones-Wassermann subfactor
associated with  ${SU}(3)_9 \subset {E_6}$ is
not commutative; \par
(3) The Kac-Wakimoto hypothesis is not true.
\endproclaim
\noindent
{\it Proof:}
(1). It follows immediately from the decompositions given above \par
and Prop.3.3 that we must have:
$$
\sigma_j \bar \rho \rho \neq \bar \rho \rho \sigma_j
$$
for some $\sigma_j$. Notice that $\sigma_j$ are irreducible subsector
of both $a_\mu$ and $\tilde a_\mu$ for some $\lambda$ by (4)
of Th.2.2. It follows from (3) of Th.2.3 
that  $\sigma_j$ commutes with all the
descendants of $a_\lambda$ and $\tilde a_\lambda$ 
for any dominant weights $\lambda$ of  ${SU}(3)$ at level 9.  Therefore
$\bar \rho \rho \in A_\rho$ must contain some irreducible sector which
don't appear  as  the
descendants of $a_\lambda$ and $\tilde a_\lambda$ 
for any dominant weights $\lambda$ of  ${SU}(3)$ at level 9. \par
(2).   Notice that
$\bar \rho \rho \bar \rho \rho
= \bar \rho \rho a_\gamma \succ  a_\gamma \succ a_{(5,1)}$, remember 
 that the dual fusion
ring is generated by the irreducible subsectors of $\bar \rho \rho$, so
$a_{(5,1)}$ appears as the subsectors in  the dual fusion
ring.  It follows
from equations (3) and (5) that  all $\sigma_j$ appear as
irreducible subsectors of the dual fusion ring. Since as shown in (1)
that 
$$
\sigma_j \bar \rho \rho \neq \bar \rho \rho \sigma_j
$$
for some $\sigma_j$, it follows that the dual fusion ring is not
commutative;\par
(3).  By using equation (4) and $\sigma_0 \sigma_2 = \sigma_1$ we get:
$$
\sigma_0 a_{(4,2)}  = 2   a_{(2,1)} + \sigma_2  +  \sigma_1   
$$. But $\bar  a_{(2,1)} =   a_{(2,1)}, \bar  \sigma_1 = \sigma_1,
\bar  \sigma_2 = \sigma_0$, it follows that the matrix
$V^{(4,2)} N_{\sigma_0}$ is not even Hermitian. By Cor.3.2, this implies
that Kac-Wakimoto hypothesis is not true in our example. \hfill Q.E.D.
%\nonident
\par
Despite the fact that (3) of Th.3.4 is the first counter-example
we find by using Cor.3.2, it is not easy to exhibit the S-matrix
explicitly.  We then find the following two examples 
 (cf. Ex. 0 and 5 in [X1]) where
the  S-matrix are much simpler and can be computed explicitly.
However, the computation of  S-matrix is not the way we find such
examples. What we did is to compute the ring structure as in the 
proof of (3) in Th.3.4 which is 
much simpler (cf. Ex. 0 and 5 in [X1]) in the next two examples, 
and by using Cor.3.2.
We shall give only the S-matrix in the next two examples. \par
{\it Example 2: (cf. Ex.0 in [X1])}
Take the conformal inclusion $SU(3)_3\subset SO(8)$. Let us label 
the weights of SU(3) as
$(\lambda_1, \lambda_2)$ with $\lambda_1 \geq \lambda_2$. The convention
is that $(0,0)$ label the vacuum or the trivial representation. 
Let $a=(2,1)$.
Let $v$
denote the vector representation of SO(8). By simple calculations one
finds:
$S_{aa} = -1/2, S_{vv} = 1/2$.
So $S_{aa}  S_{vv} = -1/4 < 0$. But it is easy to see that 
$b_{av} = 1$ .
\par
{\it Example 3: (cf. Ex.5 in [X1])}
Take the conformal inclusion $SU(4)_2 \subset SU(6)$.
The level 1 weights of $SU(6)$ are in one-to-one correspondence
with ${\Bbb Z_6} = \{\omega^i, i=0,1,...5 \}$. Let $\mu = (1,1,0)$.
One checks easily that $S_{\mu\mu} = \frac{1}{\sqrt 6}, 
S_{\omega \omega} = \frac{1}{\sqrt 6} \exp (\frac{2\pi i}{6})$.
It follows that 
$S_{\mu\mu} \bar S_{\omega \omega} = \frac{1}{ 6}  \exp (\frac{-2\pi i}{6})$.
But $b_{\omega \mu} =1$. So this is another counter example. \par
The example in Th.3.4 appears in [Kac3] where the authors claim that
the Kac-Wakimoto hypothesis has been checked to be true.  Our theorem
shows that this is not the case. However, this doesn't affect the main
purposes of  [Kac3] which is the computation of the branching
coefficients.  What we have proved is that the algorithm in  [Kac3]
is not true in general.  But in specific examples, including all three
counter examples above, the branching coefficients are easily determined
by considerations of [KW].\par 
\heading \S 4. Level-Rank duality \endheading
Level-Rank duality has been explained by different methods in 
[GW], [Sa] and [Ts].  The one that is close in spirit to our approach
is [Ts].
\subheading{4.1.  A conformal inclusion}
We shall be interested in the following  conformal inclusion:
$$
L(SU(m)_n \times SU(n)_m)  \subset \ L \ SU(nm)
$$.
In the classification of conformal inclusions in [GNO], the
above conformal inclusion corresponds to the Grassmanian
$SU(m+n)/SU(n)\times SU(m)\times U(1)$.\par
Let $\Lambda_0$ be the vacuum representation of $LSU(nm)$ on Hilbert space
$H^0$.  The decomposition of $\Lambda_0$ under $L(SU(m) \times SU(n))$ is
known, see, e.g. [Itz].  To describe such a decomposition, let us
prepare some notation.  We shall use $\dot S$ to denote the
$S$-matrices of $SU(m)$, and $\ddot S$ to denote the $S$-matrices
of $SU(n)$.  The level $n$ (resp. $m$) weight of $LSU(m)$ (resp. $LSU(n)$)
will be denoted by $\dot \lambda$ (resp. $\ddot \lambda$). \par
We start by
describing $\dot P_+^n$ (resp. $\ddot P_+^m$), i.e. the highest weights
of level $n$ of $LSU(m)$ (resp. level $m$ of $LSU(n)$).
  
$\dot P_+^n$ is the set of weights
$$
\dot \lambda = \tilde k_0 \dot \Lambda_0 + \tilde k_1 \dot \Lambda_1 +
\cdots + \tilde k_{m-1} \dot \Lambda_{m-1}
$$
where $\tilde k_i$ are non-negative integers such that
$$
\sum_{i=0}^{m-1} \tilde k_i = n
$$
and $\dot \Lambda_i = \dot \Lambda_0 + \dot \omega_i$, $1 \leq i \leq m-1$,
where $\dot \omega_i$ are the fundamental weights of $SU(m)$.
 
Instead of $\dot \lambda$ it will be more convenient to use
$$
\dot \lambda + \dot \rho = \sum_{i=0}^{m-1} k_i \dot \Lambda_i
$$
with $k_i = \tilde k_i + 1$ and $\overset m-1 \to{\underset i=0 \to \sum}
k_i = m + n$.  Due to the cyclic symmetry of the extended Dykin diagram
of $SU(m)$, the group $\Bbb Z_m$ acts on $\dot P_+^n$ by
$$
\dot \Lambda_i \rightarrow \dot \Lambda_{(i+ \dot \mu)\mod m}, \quad
\dot \mu \in \Bbb Z_m.
$$
Let $\Omega_{m,n} = \dot P_+^n / \Bbb Z_m$.  Then there is a natural
bijection between $\Omega_{m,n}$ and $\Omega_{n,m}$ (see \S2 of
[Itz]).  The proof given in [Itz] is very clear and let us 
repeat it here since it will be important later on.  The idea is to
draw a circle and divide it into $m+n$ arcs of equal length.  To
each partition $\sum_{0\leq i\leq m-1} k_i = m+n$ there corresponds
a "slicing of the pie" into $m$ successive parts with 
angles $2\pi k_i/(m+n)$, drawn with solid lines. We choose this
slicing to be clockwise.  The complementary slicing in broken lines 
(The lines which are not solid) defines a partition of $m+n$ into
$n$ successive parts, $\sum_{0\leq i\leq n-1} l_i = m+n$. We choose
the later slicing to be counterclockwise, and it is easy to see that
such a slicing corresponds uniquely to an element of  $\Omega_{n,m}$. 
\par 
 
We shall parameterize the bijection by a map
$$
\beta : \dot P_+^n \rightarrow \ddot P_+^m
$$
as follows.  Set
$$
r_j = \sum^m_{i=j} k_i, \quad 1 \leq j \leq m
$$
where $k_m \equiv k_0$.  The sequence $(r_1, \ldots , r_m)$ is decreasing,
$m + n = r_1 > r_2 > \cdots > r_m \geq 1$.  Take the complementary
sequence $(\bar r_1, \bar r_2, \ldots , \bar r_n)$ in $\{ 1, 2, \ldots ,
m+n \}$ with $\bar r_1 > \bar r_2 > \cdots > \bar r_n$.  Put
$$
S_j = m + n + \bar r_n - \bar r_{n-j+1}, \quad 1 \leq j \leq n.
$$
Then $m + n = s_1 > s_2 > \cdots > s_n \geq 1$.  The map $\beta$ is
defined by
$$
(r_1, \ldots , r_m) \rightarrow (s_1, \ldots , s_n).
$$
The following lemma summarizes what we will use .  For the proof,
see Th.1 of [Itz].
 
\proclaim{Lemma 4.1} 
%\rm{(1)}  
Let $\dot Q$ be the root lattice of
$SU(m), \
\dot \Lambda_i, \ 0 \leq i \leq m-1$ its fundamental weights and  
$\dot Q_i = (\dot Q + \dot \Lambda_i) \cap \dot P_+^n$.  
Let $\Lambda \in {\Bbb Z_{mn}}$ denote a level 1 highest weight
of $SU(mn)$ and $\dot \lambda \in \dot Q_{\Lambda \text{\rm mod} m}$.
 Then there exists a unique $\ddot \lambda
\in \ddot P_+^m$ with $\ddot \lambda = \mu \beta(\dot \lambda)$
for some unique $\mu \in \Bbb Z_n$ such that $H_{\dot \lambda} \otimes
H_{\ddot \lambda}$ appears once and only once in $H^\Lambda$.
The map $\dot \lambda \rightarrow \ddot \lambda = \mu \beta(\dot \lambda)$
is one-to-one.
Moreover, $H^\Lambda$, as representations
of $L(SU(m) \times SU(n))$, is a direct sum of all such $H_{\dot \lambda}
\otimes H_{\ddot \lambda}$.
\endproclaim
\demo{Proof}
By Th.1 of [Itz],
only the fact that the  
map $\dot \lambda \rightarrow \ddot \lambda = \mu \beta(\dot \lambda)$
is one-to-one needs to be proved.  It follows by the proof of Th.1
in [Itz] and the complete symmetry between $\dot \lambda$ and $\ddot \lambda$
that Th.1 of [Itz] remains true with  $\dot \lambda$ and  $\ddot \lambda$
exchanged. This implies the bijection. \hfill Q.E.D.
\enddemo
We shall denote by $ b(\Lambda, \dot \lambda \otimes \ddot \lambda)$
the multiplicity of  $H_{\dot \lambda}
\otimes H_{\ddot \lambda}$ in  $H^\Lambda$ and $\Lambda_0$ the
vacuum representation. 
Let $\rho$ be the sector such that
$\rho(M) \subset M$ (M is the unique hyperfinite $III_1$ factor) is the
Jones-Wassermann subfactor associated with
the conformal inclusion $
L(SU(m)_n \times SU(n)_m)  \subset \ L \ SU(nm)$.  Then
 $\rho \bar \rho = \gamma_{\Lambda_0}
:= \sum_{\dot \lambda} b(\Lambda_0, \dot \lambda \otimes \ddot \lambda)
\dot \lambda \otimes \ddot \lambda$.  Here 
$\dot \lambda \otimes \ddot \lambda$ is the sector corresponding to 
the representation   $H_{\dot \lambda}
\otimes H_{\ddot \lambda}$ of  $L(SU(m) \times SU(n))$. As before, we 
shall reserve 1 to denote the identity or the vacuum sector.  We will
be interested in the ring homomorphism $\dot \lambda \otimes \ddot \lambda
\rightarrow a_{\dot \lambda \otimes \ddot \lambda}$.
The proof of Theorem 2.2, 2.3 applies to the present case without 
modifications and we have:
\proclaim{Theorem 4.1}
(1). The ring homomorphisms 
 $\dot \lambda \rightarrow  a_{\dot \lambda \otimes 1},
\ddot \lambda \rightarrow  a_{1\otimes \ddot \lambda  }$ are
embeddings; \par
 (2). We have $ a_{\dot \lambda \otimes 1} = \sigma_\Lambda 
 a_{1\otimes \bar {\ddot \lambda}}$ where $\Lambda$ and 
$\ddot \lambda = \mu \beta (\dot \lambda)$ are as in Lemma 4.1.
\endproclaim
\demo{Proof}
(1). 
It is sufficient to show that $ a_{\dot \lambda \otimes 1},
 a_{1\otimes \ddot \lambda}$ are irreducible sectors.  
 By (2) of Th.2.2 and Frobenius duality we have:
$$
\align
\langle  a_{\dot \lambda \otimes 1},  a_{\dot \lambda \otimes 1} \rangle
& \leq  \langle  a_{\dot \lambda \otimes 1} \bar \rho, 
 a_{\dot \lambda \otimes 1} \bar \rho \rangle \\
& =  \langle \bar \rho \dot \lambda \otimes 1, 
\bar \rho \dot \lambda \otimes 1 \rangle \\
& =  \langle \rho \bar \rho,  \dot \lambda \bar {\dot \lambda} \otimes 1   
\rangle = 1
\endalign
$$
where in the last step we have used the fact that $b(\Lambda_0, \dot 
\lambda, 1)=1$ iff $\dot \lambda =1$ which follows from Lemma 4.1. 
The proof that $ a_{1\otimes \ddot \lambda}$ is irreducible is similar.
\par
(2). By a similar proof as in (1) we have that both 
$ a_{\dot \lambda \otimes 1}$ and $ \sigma_\Lambda
 a_{1\otimes \bar {\ddot \lambda}}$ are irreducible sectors. So to prove (2)
it is sufficient to show $\langle  a_{\dot \lambda \otimes 1},
\sigma_\Lambda
 a_{1\otimes \bar {\ddot \lambda}} \rangle =1$.   By (3), (4) of
Th.2.2 and Frobenius duality we have:
$$
\align
\langle  a_{\dot \lambda \otimes 1},
\sigma_\Lambda
 a_{1\otimes \bar {\ddot \lambda}} \rangle  
& = \langle  a_{\dot \lambda \otimes 1}  a_{1\otimes  \ddot \lambda},
\sigma_\Lambda \rangle \\
& = \langle  a_{\dot \lambda \otimes \ddot \lambda}, \sigma_\Lambda \rangle \\ 
& = \langle \dot \lambda  \otimes  \ddot \lambda, \gamma_\Lambda \rangle =1
\endalign
$$
where in the last step we have used Lemma 4.1. 
%\enddemo
\hfill. Q.E.D. \par
\enddemo
For another application of the conformal inclusions considered in
this section, see [X4].
\subheading {4.2.  Two series of subfactors }
Let $\dot \rho$ (resp.$\ddot \rho$) be the sectors corresponding to  the
Jones-Wassermann subfactors associated with the conformal inclusions
$SU(n+2)_n \subset SU((n+1)(n+2)/2)$ (resp. 
$SU(n)_{n+2} \subset SU((n+1)n/2)$). It is observed in [X3] that in the
case $n=3$, the two subfactors are closely related.  The aim of this
section is to show that this is true in general. \par
We shall be using the setup of \S4.1 with $m=n+2$. 
Let $\dot \rho \bar {\dot \rho} = \dot \gamma $ (resp. 
 $\ddot \rho \bar {\ddot \rho} = \dot \gamma $).  Then
$ \dot \gamma = \sum_{\dot \lambda} b(\Lambda_0'', \dot \lambda)  \dot \lambda
$ (resp. $ \ddot \gamma 
= \sum_{\ddot \lambda} b(\Lambda_0', \ddot \lambda)  \ddot \lambda
$ ), where $ b(\Lambda_0'', \dot \lambda)$ 
(resp. $ b(\Lambda_0', \ddot \lambda)$)
are branching coefficients and $\Lambda_0''$ (resp. $\Lambda_0'$) is the
vacuum representation of $L SU((n+1)(n+2)/2)$ (resp. $L SU((n+1)n/2)$).
Fortunately the branching coefficients are worked out in [LL].
The result is the following (cf. Th.1.2 and Th.2.1 of [LL]):
Define the group $({\Bbb Z_2})^m_{\text{\rm even}} := \{(s_{11},...
s_{1m} )\in ({\Bbb Z_2})^m: \sharp \{ i, s_{1i} = -1\} \in 2 {\Bbb N} \}$ .
Define $c(s_1) := \sum_{i,s_{1i} =1} (m-i)$.  Let $\mu_1$ (resp. $\mu_2$)
be the generator of the center of $SU(m)$ (resp. $SU(n)$).
Let $\{ a_1,...a_m \} = \{s_{1i}(m-i) , i=1,...m \}$  with $ a_1,...a_m$
in decreasing order.
Define a weight 
\footnotemark\footnotetext{ This is to be compared to the similar
expressions in \S 2 of [LL].  The weights in this paper differ from
the weights of [LL] by the Weyl vector according to our convention in
\S4.1} $\sigma_{s_1} (\rho_1) := (2m-2 + a_m -a_1) \dot \Lambda_0
+ \sum_{1\leq i \leq m-1} (a_i-a_{i+1}) \dot \gamma_i$ .  
Similarly define $ s_2 \in ({\Bbb Z_2})^n$ to be $s_2=  (s_{21},...
s_{2n})$. Let  $\{ b_1,...b_n \} = \{s_{2i}(n+1-i) , i=1,...n \}$  
with $ b_1,...b_n$
in decreasing order.
Define a weight $\sigma_{s_2} (\rho_1) := (2n+2 + b_n -b_1) \ddot \Lambda_0
+ \sum_{1\leq i \leq n-1} (b_i-b_{i+1}) \ddot \gamma_i$. \par

Then it follows from Th.2.1  [LL] that 
 $b(\Lambda_0'', \dot \lambda)=1$ iff 
$$
\dot \lambda = \mu_1^{2k_1} (\sigma_{s_1} (\rho_1))
$$, with $ 0\equiv c(s_1) + 2k_1(m-1) \text{\rm mod} \ \ m(m-1)/2$.
Notice in the notaton of Lemma 4.1, such a $\dot \lambda $ belongs
to $\dot Q_0$, and there exists a unique $\ddot \lambda = \mu \beta (
 \dot \lambda)$ such that $b(\Lambda_0, \dot \lambda, \ddot \lambda) =1$.
We claim that $b(\Lambda_0', \ddot \lambda) =1$.  This is completed in
two steps: \par
First we show that $\beta(\sigma_{s_1} (\rho_1)) = \mu_2^k (\sigma_{s_2}
(s_2(\rho_2))$ for some $k$ and $s_2$.  Notice $0$ appears in the sequence
 $\{ a_1,...a_m \}$ and either $a_1=m-1$ or $a_m = -(m-1)$ but one can not
have both $a_1=m-1$ and  $a_m = -(m-1)$ . 
If  $a_1=m-1$, we define  $\{ b_1',...b_n' \}$ to be the complement
sequence of  $\{ a_1,...a_m \}$ in $0,1,... m-1, -1,-2,... -(m-2)$, 
in decreasing order.
If  $a_m = -(m-1)$,  we define  $\{ b_1',...b_n' \}$ to be the complement
sequence of  $\{ a_1,...a_m \}$ in $0,1,... m-2, -1,-2,... -(m-1)$
in decreasing order.
Define  $\{ b_1,...b_n \}$ to be $\{ -b_n',... - b_1' \}$, in decreasing
order. It is then easy to see that there is a unique $s_2 \in ({\Bbb Z_2})^n$
such that  $\{ b_1,...b_n \} = \{ s_{2i}(n+1-i), i=1,2,... n\}$.  \par
Now consider the unit circle on the $xy$-plane.  Divide it into $m+n$ arcs
of equal length.  Draw the solid lines labeled by $a_i$'s such that 
the angle between  the solid line labeled by $a_i$ and the positive $x$-axis
is $2\pi i a_i/(m+n)$.  Draw the broken lines  labeled by $b_i'$'s such that   
the angle between  the solid line labeled by $b_i'$ and the positive $x$-axis
is $2\pi i b_i'/(m+n)$.  It is easy to see that the pie formed by the 
broken lines, in anti-clockwise direction, corresponds the orbit of the
weight $ (2n+2 + b_n -b_1) \ddot \Lambda_0
+ \sum_{1\leq i \leq n-1} (b_i-b_{i+1}) \ddot \gamma_i$, which is precisely
 $\sigma_{s_2} (\rho_1)$.  By
the "slicing pie" argument in \S4.1,  it is then easy to see that 
there exists a $k$ such that 
$\beta(\sigma_{s_1} (\rho_1)) = \mu_2^k (\sigma_{s_2}
(s_2(\rho_2))$. \par
It follow that $\mu \beta (\dot \lambda) = \ddot \lambda = \mu^{k_2}
( (\sigma_{s_2}
(s_2(\rho_2))$ for some $k_2 \in {\Bbb Z}$.  So we have (cf. Page 9 of
[LL]):
$$
b(  \Lambda', \ddot \lambda) = 1
$$ for some weights $  \Lambda'$ of $SU(n(n+1))/2$.  Let us show 
that $  \Lambda'$ is the vacuum representation $  \Lambda_0'$.
Since $b(  \Lambda',\ddot \lambda) = 1, b(\Lambda_0'', \dot \lambda)=1,
b(\Lambda_0, \dot \lambda, \ddot \lambda) = 1$, it follow from [KW] that:
$h_{  \Lambda'} - h_{\ddot \lambda} \in {\Bbb Z}$,
$ h_{\ddot \lambda} +  h_{\dot \lambda} \in {\Bbb Z}$, and
$ h_{\dot \lambda} \in {\Bbb Z}$, where $h_a$ is the conformal
anomaly of the weight $a$.  So we have:
$$
h_{ \Lambda'}  \in {\Bbb Z}
$$, and it follows from (1.19) of [LL] that 
 $ \Lambda'$ is the vacuum representation $  \Lambda_0'$. \par
To summarize, we have shown that if $b(\Lambda_0'', \dot \lambda)=1$
and $\ddot \lambda = \mu \beta (\dot \lambda)$ as in Lemma 4.1, then
 $b(\Lambda_0', \dot \lambda)=1$.  Notice $ d_{\dot \lambda} =
 d_{a_{\dot \lambda \otimes 1}} =  d_{a_{1\otimes \ddot \lambda }} =
 d_{\ddot \lambda}$ by Th.4.1,  
where $d_a$ is the statistical dimension of sector $a$.
Since $\dot \lambda \rightarrow \mu \beta
  (\dot \lambda)$ is one to one, we have proved the following:
$$
\align
\sum_{\dot \lambda} b(\Lambda_0'', \dot \lambda) d_{\dot \lambda}
& = \sum_{\dot \lambda} b(\Lambda_0', \mu \beta (\dot \lambda)) 
d_{\mu \beta (\dot \lambda)} \\
& \leq   \sum_{\ddot \lambda} b(\Lambda_0', \ddot \lambda) d_{\ddot \lambda}
\endalign
$$
.  Now if we start with $\ddot \lambda$ with  $b(\Lambda_0', \ddot \lambda)=1$
and go through the previous arguments, we obtain the reverse inequality
above, therefore proving that the above inequality is actually an
equality. So we have shown  $\dot \lambda \rightarrow \mu \beta
  (\dot \lambda)$ is a one-to-one and onto map between the 
irreducible subsectors of $\dot \rho \overline{\dot \rho}$ and
 $\ddot \rho \overline{\ddot \rho}$.  Since $b(\Lambda_0', \dot \lambda) =
b(\Lambda_0', \overline{\dot \lambda)}$, the map $\phi: 
\dot \lambda \rightarrow \overline{\mu \beta
  (\dot \lambda)}$ is a one-to-one and onto map between the
irreducible subsectors of $\dot \rho \bar {\dot \rho}$ and
 $\ddot \rho \bar {\ddot \rho}$.  By (2) of Th.4.1, 
$a_{\dot \lambda \otimes 1} = a_{1\otimes \overline{\mu \beta
  (\dot \lambda)}}$.  It follows from (1) of Th.4.1 that  
the map $\phi$ is really a ring isomorphism.  Since the fusion
ring of  $\dot \rho $ (resp.  $\ddot \rho$) is generated by 
irreducible subsectors of  $\dot \rho \bar {\dot \rho}$ (resp.
 $\ddot \rho \bar {\ddot \rho}$), we have proved the following theorem:
\proclaim{Theorem 4.2}
The fusion ring of the Jones-Wassermann subfactors associated with
 the conformal inclusions
$SU(n+2)_n \subset SU((n+1)(n+2)/2)$ and 
$SU(n)_{n+2} \subset SU((n+1)n/2)$ are canonically isomorphic via
$\phi$ defined above.
\endproclaim
% \enddemo
\heading \S 5. Conclusions and questions \endheading
In this paper we have given applications of the general theory
developed in [X1].  The following questions arise naturally
from our approach: \par
(1) In all the counter examples we have about Kac-Wakimoto
hypothesis, there are always multiplicities, i.e., there is
a $\mu$ such that $b_{i\mu}\neq 0, b_{j\mu}\neq 0$ for some $i\neq j$.
It is not clear to us if there are examples in the multiplicity-free case;
\par
(2) Let $s(\mu,i)$ denote the multiplicity of $(\mu.i)$ in 
(Exp), is it true that $s(\mu,i) = b_{\mu,i}$?  If this is true,
it will imply that the ring $C_\rho$ in the conformal inclusions
considered in \S4.2 is commutative. This question seems to be related
to the $M$-algebra in [PZ];  \par
(3) In \S 4.2 we have studied the fusion ring of the Jones-Wassermann
subfactors associated with two series of conformal inclusions. It will 
be interesting to see if the ring $C_\rho$ of these  conformal inclusions
is related in simple way; \par
(4) The  conformal inclusions
considered in \S4.2 is related to the compact Hermitian symmetric spaces
$Sp(n)/U(n)$ and $SO(2n)/U(n)$.  It will be interesting to
see if the ring $C_\rho$ is related to the quantum cohomology ring
of those  Hermitian symmetric spaces. 

\heading References \endheading
\roster
\item"{[L1]}"  R. Longo, {\it Minimal index and braided subfactors}, J.
Funct. Anal., {\bf 109}, 98-112 (1992).
\item"{[L2]}"  R. Longo, {\it Duality for Hopf algebras and for subfactors},
I, Comm. Math. Phys., {\bf 159}, 133-150 (1994).
\item"{[L3]}"  R. Longo, {\it Index of subfactors and statistics of
quantum fields}, I, Comm. Math. Phys., {\bf 126}, 217-247 (1989.
\item"{[L4]}"  R. Longo, {\it Index of subfactors and statistics of 
quantum fields}, II, Comm. Math. Phys., {\bf 130}, 285-309 (1990).
\item"{[L5]}"  R. Longo, Proceedings of International Congress of
Mathematicians, 1281-1291 (1994).
\item"{[LR]}"  R. Longo and K.-H. Rehren, {\it Nets of subfactors},
Rev. Math. Phys., {\bf 7}, 567-597 (1995).
\item"{[Po]}" S.Popa, {\it Classification of subfactors and of their
endomorphisms}, CBMS Lecture Notes Series, 86.
\item"{[PP]}" M.Pimsner and S.Popa, 
{\it Entropy and index for subfactors}, \par
Ann. \'{E}c.Norm.Sup. {\bf 19},
57-106 (1986).
\item"{[W1]}"  A. Wassermann, {\it Operator algebras and Conformal
field theories III},  to appear.
\item"{[W2]}"  A. Wassermann, Proceedings of International Congress of 
Mathematicians, 966-979 (1994).
\item"{[PZ]}"  V. B. Petkova and J.-B. Zuber, {\it From CFT to graphs},
hep-th-9510198.
\item"{[Fran]}"  P. Di Francesco and J.-B. Zuber, {\it Integrable lattice
models associated with $SU(N)$}, Nucl. Phys. B, {\bf 338}, 602-623 (1990).
\item"{[Ka1]}"  Y. Kawahigashi, {\it Classification of paragroup actions
on subfactors}, Publ. RIMS, Kyoto Univ., {\bf 31} 481-517 (1995). 
\item"{[GHJ]}"  F. M. Goodman, P. de la Harpe and V. Jones, {\it Towers of
algebras and Coxeter graphs}, MSRI publication, no. 14.
\item"{[PS]}"  A. Pressley and G. Segal, {\it Loop groups}, Clarendon Press,
1986.
\item"{[GNO]}"  P. Goddard, W. Nahm and D. Olive, {\it Symmetric spaces,
Sugawara's energy momentum tensor in two dimensions and free fermions},
Phy. Lett., {\bf 160B}, 111-116 (1985)
\item"{[KW]}"  V. G. Kac and M. Wakimoto, {\it Modular and conformal
invariance constraints in representation theory of affine algebras},
Advances in Math., {\bf 70}, 156-234 (1988).
\item"{[Kac2]}"  V. G. Kac, {\it Infinite dimensional algebras}, 3rd Edition,
Cambridge University Press, 1990.
\item"{[Kac3]}"  V. G. Kac and M. Sanielevici, {\it Decompositions of
representations of exceptional affine algebras with respect to
conformal subalgebras},
Phys.Rev.D 37, {\bf 70}, 2231-2237 (1988).
\item"{[X1]}"  F. Xu, {\it New braided endomorphisms from conformal
inclusions},  \par
Comm. Math. Phys., (1997) (in press)
\item"{[X2]}"  F. Xu, {\it Generalized Goodman-Harper-Jones construction
of subfactors, I},  Comm. Math. Phys., {\bf 184}, 475-491 (1997).
\item"{[X3]}"  F. Xu, {\it Generalized Goodman-Harpe-Jones
construction of subfactors}, II,  Comm. Math. Phys., {\bf 184}, 493-508 (1997).
\item"{[X4]}" F.Xu, {\it Jones-Wassermann subfactors for
Disconnected Intervals}, \par q-alg 9704003.
\item"{[Y]}"  S. Yamagami, {\it A note on Ocneanu's approach to Jones
index theory}, Internat. J. Math., {\bf 4}, 859-871 (1993).
\item"{[GW]}"  F. M. Goodman and H. Wenzl, {\it Littlewood Richardson
coefficients for Hecke algebras at roots of unity}, Adv. Math., 1990.
\item"{[S]}"  H. Saleur, {\it Level-rank duality}, Nucl. Phys. B.,
{\bf 363}, 177-192 (1992).
\item"{[Ts]}"  T.Nakanishi and A.Tsuchiya, 
{\it Level-rank duality of WZW models in conformal field
theory},  Comm. Math. Phys.,
{\bf 144}, 351-372 (1992).
\item"{[DB]}"  D.Bisch,
{\it On the structure of finite depth subfactors},  Algebraic
Methods in Operator Theory, Birkh\"{a}user, 175-194.
\endroster
\enddocument